# Opportunities and Challenges Applying Functional Data Analysis to the Study of Open Source Software Evolution

**Katherine J. Stewart, David P. Darcy and Sherae L. Daniel**

*Abstract.* This paper explores the application of functional data analysis (FDA) as a means to study the dynamics of software evolution in the open source context. Several challenges in analyzing the data from software projects are discussed, an approach to overcoming those challenges is described, and preliminary results from the analysis of a sample of open source software (OSS) projects are provided. The results demonstrate the utility of FDA for uncovering and categorizing multiple distinct patterns of evolution in the complexity of OSS projects. These results are promising in that they demonstrate some patterns in which the complexity of software decreased as the software grew in size, a particularly novel result. The paper reports preliminary explorations of factors that may be associated with decreasing complexity patterns in these projects. The paper concludes by describing several next steps for this research project as well as some questions for which more sophisticated analytical techniques may be needed.

## 1. INTRODUCTION

Software complexity is a crucial factor in many important outcomes of the software development process, including defect rates, maintainability, security, and reliability (Kemerer, 1995). As these outcomes are often viewed as being at the root of what makes for "better" or "worse" software, complexity has been seen as a key contributor to overall software quality (Prahalad and Krishnan, 1999). Understanding software complexity and how it may be managed throughout the software lifecycle is therefore of great interest to software developers and researchers.

Much of the thinking on the evolution of software has come from Lehman's laws of software evolution (Belady and Lehman, 1976). Relevant insights include that software will tend to grow in size and functionality from release to release in order to remain valuable to users. Similarly, complexity of the software will tend to increase through the life of a project, unless it is actively managed. Much of the empirical work based on the laws has largely supported them (Kemerer and Slaughter, 1999), although that work has mostly taken place in the closed source context, and has often been limited to studying a single system or a small set of systems.

*Katherine J. Stewart is Assistant Professor of Information Systems, Department of Decision and Information Technologies, Robert H. Smith School of Business, University of Maryland, College Park, Maryland 20742, USA e-mail: kstewart@rhsmith.umd.edu. David P. Darcy is Assistant Professor of Information Systems, Department of Decision and Information Technologies, Robert H. Smith School of Business, University of Maryland, College Park, Maryland 20742, USA e-mail: ddarcy@rhsmith.umd.edu. Sherae L. Daniel is a Doctoral Candidate, Department of Decision and Information Technologies, Robert H. Smith School of Business, University of Maryland, College Park, Maryland 20742, USA e-mail: sdaniel@rhsmith.umd.edu.*

While most work to date has examined software complexity in the context of closed source software development, managing complexity could be at least as important in the open source development context. OSS is software released under a license approved by the Open Source Initiative (OSI, see www.opensource.org). The main OSI licensing requirement is that source code be available when the software is distributed. As an example of an OSS project, see the statistical application R that we used in our analyses (www.r-project.org). Because source code is available, OSS is open to a wide audience for use, inspection, contributions and modifications. Software complexity may be especially important in this context because of the potentially higher fluidity in membership of the project team (i.e., developers may join and leave more freely than in most closed development contexts and new members will be able to more quickly contribute if complexity is minimized), the lack of or lag in formal design specifications and documentation to aid in developing, understanding and maintaining the software (Scacchi, 2002) and the support concerns that have been voiced by many potential OSS adopters (Smith, 2002).

The goal of this research is to uncover and characterize patterns of evolution in the complexity of OSS projects. This is a first step in investigating the antecedents and consequences of different evolutionary patterns, which we hope may provide insights for improved project management in both open and closedsource software project contexts. The objective of this paper is to describe how functional data analysis may be applied toward achieving the research goals. Given this objective, our discussion of the large literature on software development and evolution is extremely limited. While we do not delve far into theoretical issues surrounding software evolution, we do discuss some of the drawbacks to using more established statistical techniques to analyze evolution. This serves as an introduction to the opportunities for uncovering patterns in evolution using functional data analysis (FDA), several challenges we encountered in analyzing software complexity data, and the approaches we considered and selected to address those challenges. The final section of the paper discusses the preliminary conclusions drawn from the FDA results, and the analytical limitations that remain. To provide the backdrop for this discussion, the next two sections provide a brief overview of literature on software complexity and a description of the data collected for the study.

## 2. THEORETICAL BACKGROUND

Several aspects of software complexity have been studied. Algorithmic complexity, for example, examines the machine resources, such as time necessary to solve a given problem. Structural software complexity has been defined as "the organization of program elements within a program" (Gorla and Ramakrishnan, 1997, page 191). Software designed to solve a particular problem can be structured in many different ways. Different structures may lead to significant variations in the amount of effort required in implementing or maintaining a software solution. For example, there is evidence to suggest that the complexity of the Linux kernel is so high that future changes will be difficult (Yu, Schach et al., 2004). While algorithmic complexity is important in terms of providing the requisite machine resources, structural complexity has an impact on the implementation and maintenance (i.e., human) resources necessary to provide software. The availability of machine resources continues to increase, and their cost to decrease, at exponential rates, while the availability and cost of human resources necessary to implement and maintain software is more constant. Because the human resources involved in building and maintaining software represent an increasing portion of the total cost of software, we chose to focus on structural complexity.

In this paper, two dimensions of structural complexity of software are considered: coupling and cohesion (Chidamber, Darcy and Kemerer, 1998). Coupling (Cpl) is the degree to which a program element is "related" to other elements; the higher the average coupling of elements in the software, the more complex it is considered to be. Coupling measures tend to use absolute scales; for each instance of coupling by a program element to another element, the measure is incremented by 1. Cohesion is the degree to which the content within a program element is related; the higher the average cohesion of elements in the software, the less complex it is considered to be. Cohesion measures tend to use a percent scale; more cohesive program elements will have values closer to 100%. For consistency, we use lack of cohesion as a measure (abbreviated as LCoh) such that increases in Cpl and LCoh both represent increases in complexity.

Coupling and cohesion have been argued to be inversely correlated such that managing one may result in a trade-off for managing the other (Chidamber, Darcy and Kemerer, 1998). Particularly when considering a program in its entirety rather than an individual program element, different design choices will often impact both coupling and cohesion, frequently reducing one at the expense of increasing the other. Therefore it is the interaction of coupling and cohesion, rather than either one alone, that may best determine software maintenance effort (Darcy, Kemerer et al., 2005). To take both measures into account we followed Darcy, Kemerer et al. (2005) and calculated a cross-term, multiplying coupling by cohesion. This cross-term is the main variable of interest in the remainder of this paper, and it is referred to hereafter as "CplXLCoh."

## 3. RESEARCH DESIGN AND DATA COLLECTION

As part of a larger project, data were collected on 105 OSS projects hosted online at Sourceforge (sf.net). Sourceforge is the largest OSS repository, currently hosting over 100,000 projects and over one million registered users. OSS projects use Sourceforge to manage development and to make releases of their software available. In order to limit the variance in structural complexity driven by factors other than software design choices, we limited our data collection to projects that use only the Java programming language and are listed in the Internet and System Networking domains. By examining only projects written for the Internet and System Networking domains, we limit variance in structural complexity that may be driven by the underlying problem type. Similarly, by exploring only Java projects, we limit variance in structural complexity that may be driven by the programming language. We further limited the data collection by only including these projects that use an OSI approved license to ensure that only software that is truly open source was included. Finally, to be included a project had to have posted at least one file on the Sourceforge site as of the time of our initial project selection, Fall 2002. This criterion was applied to screen out projects that had not produced any code because projects may be listed on Sourceforge soon after their inception, before any software has actually been produced. Data were collected on the published release history of each project that met the screening criteria. Each release of each project was analyzed to calculate CplXLCoh. The size of each release was measured using a calculation of the number of lines of code (LOC).

Figure 1 shows the CplXLCoh measure for each release for three of the projects in the sample, and Table 1 provides the data associated with one of these projects (the squares in Figure 1). The $x$-axis represents calendar time; it starts at January 2000 because that is the earliest release for any project in the sample. Figure 1 is included to graphically depict several challenges that had to be overcome to analyze the data. These are that (1) the projects have different starting points, (2) the projects have different ending points, (3) the project histories are of different lengths, (4) the projects have different numbers of releases (i.e., data points) and (5) the projects span different levels of complexity. The next section discusses in greater detail how each of these issues was accommodated to derive functional objects for each project.

## 4. CHALLENGES AND STRATEGIES FOR EXAMINING SOFTWARE EVOLUTION

Like the online auction setting described by Jank and Shmueli (2005), studying software evolution poses several challenges in terms of both the nature of the data available and the insights we seek to generate using it. The data have all of the unique aspects of auction bid data outlined by Jank and Shmueli (2005), including that data from releases form a time series, but they are unequally spaced and occur with differing frequencies over varying amounts of time for different projects. An additional challenge is posed because unlike auctions where there is a final winning price, OSS software projects generally have no defined end point.

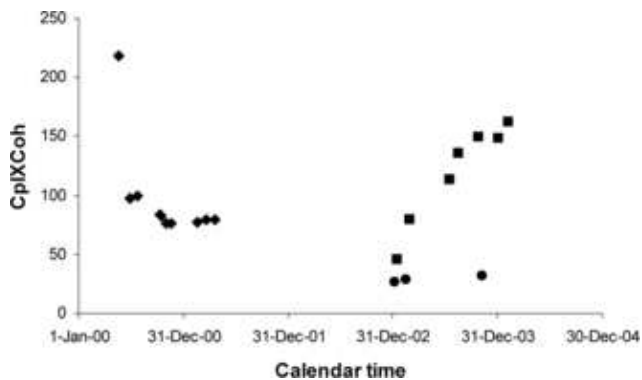

FIG. 1. *Three sample project release histories.*

TABLE 1
*Sample project data*

| Project ID | Release date | LOC | CplXLCoh |
|------------|--------------|------|----------|
| 3064 | 17-Jan-03 | 4,901 | 45.71 |
| 3064 | 2-Mar-03 | 5,449 | 79.31 |
| 3064 | 16-Jul-03 | 6,775 | 113.83 |
| 3064 | 16-Aug-03 | 10,915 | 135.98 |
| 3064 | 25-Oct-03 | 13,516 | 149.15 |
| 3064 | 4-Jan-04 | 13,991 | 148.65 |
| 3064 | 7-Feb-04 | 14,892 | 162.30 |

While most empirical work on software evolution has focused on studying releases of a single system over time (Kemerer and Slaughter, 1999), we wish to analyze multiple sets of time-series data on OSS projects in order to gain insights that may be applicable to other projects as opposed to analyzing a single time series in order to gain insight into a single project. Thus we turned to functional data analysis, which allows for capturing the dynamics of the projects over time, representing those dynamics using polynomial pieces to create a function for each project, and then using these functions to analyze patterns in the projects.

**Project Starting Point**

The fact that projects in the sample begin on different calendar days presents a question concerning how to align the data so that functions derived for different projects are comparable. As implied in Figure 1, the projects in the sample do not present an obvious, "natural" alignment. Aligning the starting points of the curves is easy if we ignore the calendar time at which projects commence and simply place the first release at the origin on the *x*-axis and plot each curve from that point forward. This approach mimics that taken in studies that have used age on the *x*-axis (Ramsay and Silverman, 2002), where the date of birth of a person can be safely ignored because the interest is in comparing changes across lifespan that are assumed independent of date of birth.

A potential difficulty in applying this approach in this context is that the first public releases for different projects do not necessarily represent similar points in the evolution of the projects. For example, some OSS projects may release an initial version very early in development when only a small portion of the software has been created whereas others may release an initial version only after most of the intended functionality has been created and tested. Thus "day 1" for OSS projects is not equivalent to "age 1" for a person. Unfortunately, however, the data do not provide any good indication of the development stage of the project at its first release; thus aligning the first release of each project as the starting point was the best approach available. This approach has implications for the nature of the conclusions that can be drawn from analysis, which will be discussed below.

**Project Ending Point**

Just as there is no reason to believe the first releases represent comparable points in the development of the software of each project, there is no basis for assuming that the last releases represent comparable points. As can be seen in Figure 1, the data span different amounts of time, and for those projects with final releases close to the last point on the *x*-axis, there is no way to know what may have occurred in the project development after the end of the period.

In order to generate curves that are as equivalent as possible in representing similar periods of development for each project, we limited the dataset to include only those releases that occurred within two years of the initial release of a project (730 days). We return to the implications of this decision in discussion of the results, below.

**Creating Comparable Project Data**

Having decided to consider data for releases during the first two years after the initial release of each project in order to align starting and ending points, there remained the issue that some projects had data spanning periods that were significantly shorter than others. In other words, though each project had one release at least two years prior to data collection, many projects have long periods of time with no releases. For example, the project represented by circles in Figure 1 had three releases spanning 301 days and then nothing since. Whereas we cannot be certain that this indicates development on the project stopped (e.g., the project could have released a new, different version one day after our data collection), this is how we chose to interpret it. We thus assumed the level of complexity was constant after the last release, for the remainder of the two-year period. This treatment of the data is consistent with viewing the software from the user perspective in that even if development was actively under way, the complexity of the software available to users was unchanged.

Having expanded data for short projects so that each project was the same length, the final issue in preparing the data to generate curves was to consider alignment of points along the curves. Projects produced different numbers of releases, making it difficult to align change points and produce comparable functional representations across projects. Following prior work (Ramsay and Silverman, 2002; Jank and Shmueli, 2005), this issue was addressed by interpolating the data to create equivalent numbers of observations. We expanded the data to create values of complexity for every project for every day of the two-year period using a step function. The values calculated from the first release of the project were assigned for every day until the second release, and then values from the second release were assigned for every day until the third release, and so forth. An alternative approach would be to assume a linear change in complexity between releases; however, we chose the former approach to maintain consistency with the manner in which the data were extended to cover the entire two-year period and because, again, from the user perspective, changes occur at the discrete points when new releases become available.

## 5. RESULTS

Because our interest is in understanding how complexity evolves as software grows over time, we screened out projects that did not have at least a 5% increase in LOC between the first and last releases in a single development stream of the project. This mainly removed projects that had only released a single version of their software, and left 59 projects for analysis. Descriptive statistics for these 59 projects are provided in Table 2.

Our analysis is divided into two main parts. The first part examines the patterns of change in actual values of complexity. Based on the findings of the first analysis, the second analysis examines the evolution of complexity values that have been standardized within projects.

### Analysis of Absolute Levels of Complexity

Figure 2 shows the same three projects as Figure 1, but with the data adjusted as described. Our next step was to use smoothing splines to obtain representations of the projects that would be more compact and reduce some of the fluctuation to allow us to discern underlying patterns. To fit smoothing splines, the overall interval of interest, in this case 2 years, is split into subintervals and a polynomial piece is fitted for each subinterval such that the polynomial pieces are fit together smoothly at the points where the subintervals meet. These points are the knots. Minimization of the penalized residual sum of squares is done in a way similar to the minimization of the least-squares operator in regression analysis (see Shmueli and Jank, 2006). Following Jank and Shmueli (2005), we sampled from the step functions in order to place knots to connect the smoothing splines across subintervals. For this analysis a parameter of 13 knots was chosen as this represents the number of project days explored (730) divided by the average number of days between releases (56). That is, on average the projects released changes every 56 days, and therefore a 56-day interval was chosen. We used the smooth.spline function to generate the curves (see ego.psych.mcgill.ca/misc/fda/software.html).

The order of the smoothing spline determines how much the curves may deviate from a flat line. Prior work employing FDA has used a cubic smoothing spline (Ramsay and Silverman, 2002), and because we have no reason to expect a higher-order function to be more appropriate, that convention is followed here. The smoothing parameter for the curves determines the degree to which the function is faithful to the observed data points. We selected a high value to generate relatively smooth curves because a primary objective of this initial study is to develop a description of broad underlying patterns, and allowing for wider fluctuations in the functions makes interpretation of the basic patterns more difficult.

Table 2
*Descriptive statistics (n = 59 OSS projects)*

| Measure | Mean | Std dev | Min | Max |
|---|---|---|---|---|
| Number of releases | 8.44 | 6.56 | 2 | 29 |
| Average release frequency (in days) | 56.08 | 52.05 | 0.00 | 238.00 |
| First release LOC | 7,031.32 | 9,589.37 | 395.00 | 52,792.00 |
| Last release LOC | 14,565.63 | 16,165.28 | 595.00 | 70,012.00 |
| First release CplXCoh | 98.28 | 54.63 | 13.06 | 239.16 |
| Last release CplXCoh | 113.64 | 58.01 | 18.36 | 360.04 |

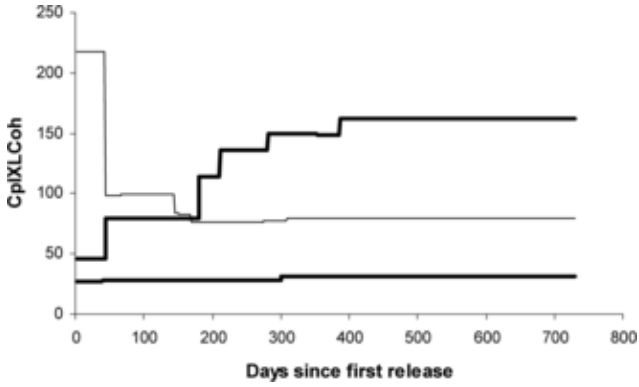

Fig. 2. *Adjusted data for the same three sample projects.*

The adjusted data for all 59 projects were used to create and plot functional objects, examine the mean functional object and explore whether there were different clusters of projects. The mean curve, shown as a solid line in Figure 3(a), is virtually flat, with only a very small positive slope. The dashed lines indicate a 95% confidence interval calculated pointwise using the fitted values. Overall, the mean curve shows very little upward movement in complexity, which is a surprising result given a large quantity of past work arguing that complexity increases as software grows (Belady and Lehman, 1976; Kemerer and Slaughter, 1999), and we specifically limited the sample to projects that increased in size (the average increase in LOC was 107%). Examining the individual project curves shows that many projects decreased in complexity over their lives, while others followed the expected pattern of increase, which sheds light on why the overall mean is flat.

Following Jank and Shmueli (2005), we used K-medoids clustering on the coefficients of the functions to separate the opposing patterns that result in an overall flat mean. K-medoids clustering is a more robust version of K-means clustering, especially with respect to outliers. The K-medoids algorithm (e.g., Kaufman and Rousseeuw, 1987; Hastie, Tibshirani and Friedman, 2001) minimizes within-cluster dissimilarity. This is done by iteratively alternating between two steps. During the first step the cluster center observations are determined based on the current data-partitions. The cluster center observation in the $k$th cluster that minimizes the total distance to the other points in the cluster is chosen.

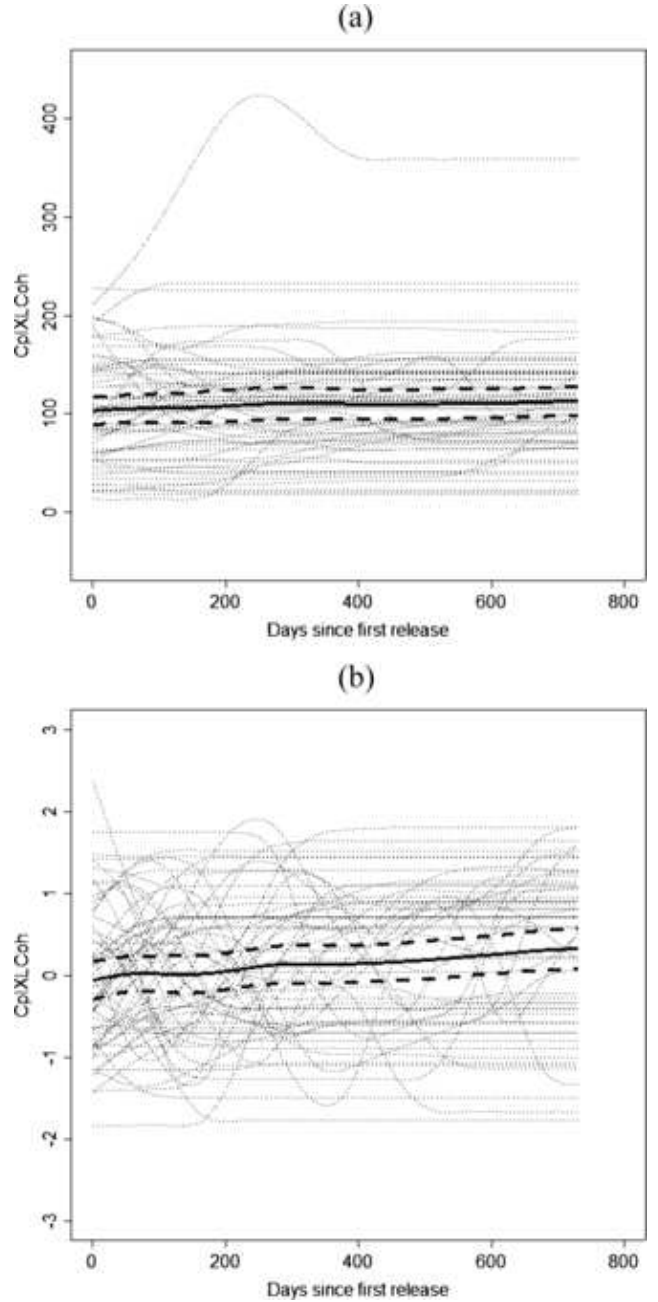

Fig. 3. *(a) Functional objects, absolute values ($n = 59$, solid mean and dashed 95% CI). (b) Functional objects, standardized values ($n = 59$, solid mean and dashed 95% CI).*

This observation is known as the medoid and is the most centrally located point in a cluster. During the second step observations are reassigned to the nearest medoid and the associated cluster. These steps continue until the assignments do not change.

When we applied cluster analysis to the data, we found that different patterns of change in evolution were overwhelmed by different absolute levels of complexity such that the results separated projects into clusters that had similar overall levels of complexity, but the mean function for each cluster was still similar in shape to the overall mean function (see Figure 4). While it is crucial to manage the absolute level of complexity of software code, it is also widely accepted that complexity is correlated with size (Chidamber, Darcy and Kemerer, 1998), and a significant correlation between LOC and CplXL-Coh (0.259, $p \leq 0.05$) in the dataset indicated consistency with this assertion. Since complexity is a function of size, which is largely determined by the scope of the problem the software is designed to address, the extent to which a project manages changes in complexity over time (e.g., minimizes increases in complexity that may result from maintenance) is of greater interest to us than the absolute level of complexity in a project. Thus we sought to adjust the data to focus on the pattern of change (i.e., the shapes of the curves) rather than the absolute level of complexity. To do this, we standardized the values of complexity for each project such that the average for every project was at 0. This allowed us a means of aligning the curves, and results of analysis on the standardized values produce more meaningful clusters in terms of having identifiably different shapes.

**Analysis of Standardized Complexity**

Parameter values for the analysis on the standardized data were the same as those used to generate Figure 3(a). The overall mean curve derived from the standardized data has the same flat shape as the curve for the absolute data [see Figure 3(b)]. However, conducting cluster analyses yields more interesting results. Because this research is exploratory in that there are no prior established patterns of complexity evolution that we expect the data to follow, we began by exploring a two-cluster solution and then increased the number of clusters until no new patterns were uncovered. The two-cluster solution [Figure 5(a)] resulted in 25 projects in cluster 1 that show an overall decline in comlexity over approximately the first 350 days, then flatten out before starting a slight upward trend toward the end of the observation period. The 34 projects in cluster 2 have an overall upward in complexity with the fastest rate of increase over roughly the first 300 days. The main point of interest in the two-cluster solution is that it results in one cluster with projects whose complexity increases and one cluster for projects whose complexity decreases. This is consistent with the overall flat mean for the entire set; however, it is a very surprising result in the sense that no prior empirical work in either the open or closed source software context has identified a development pattern in which complexity decreases as size is increased.

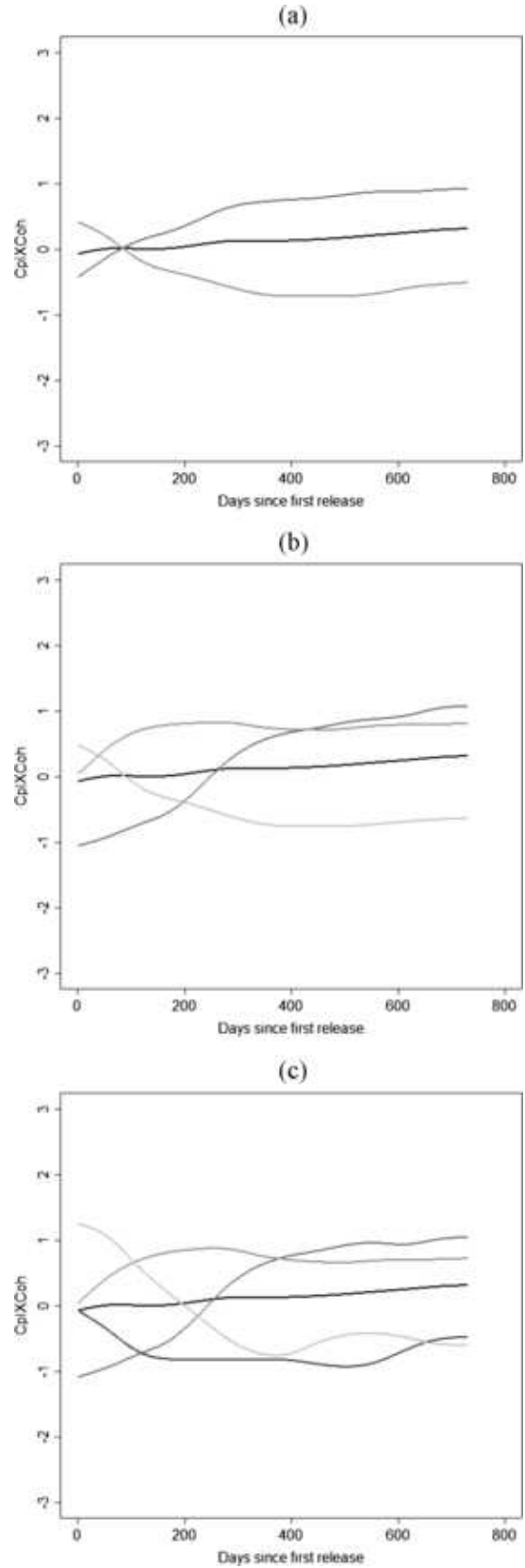

FIG. 5. (a) Two-cluster mean functions (solid overall mean). (b) Three-cluster mean functions (solid overall mean). (c) Four-cluster mean functions (solid overall mean).

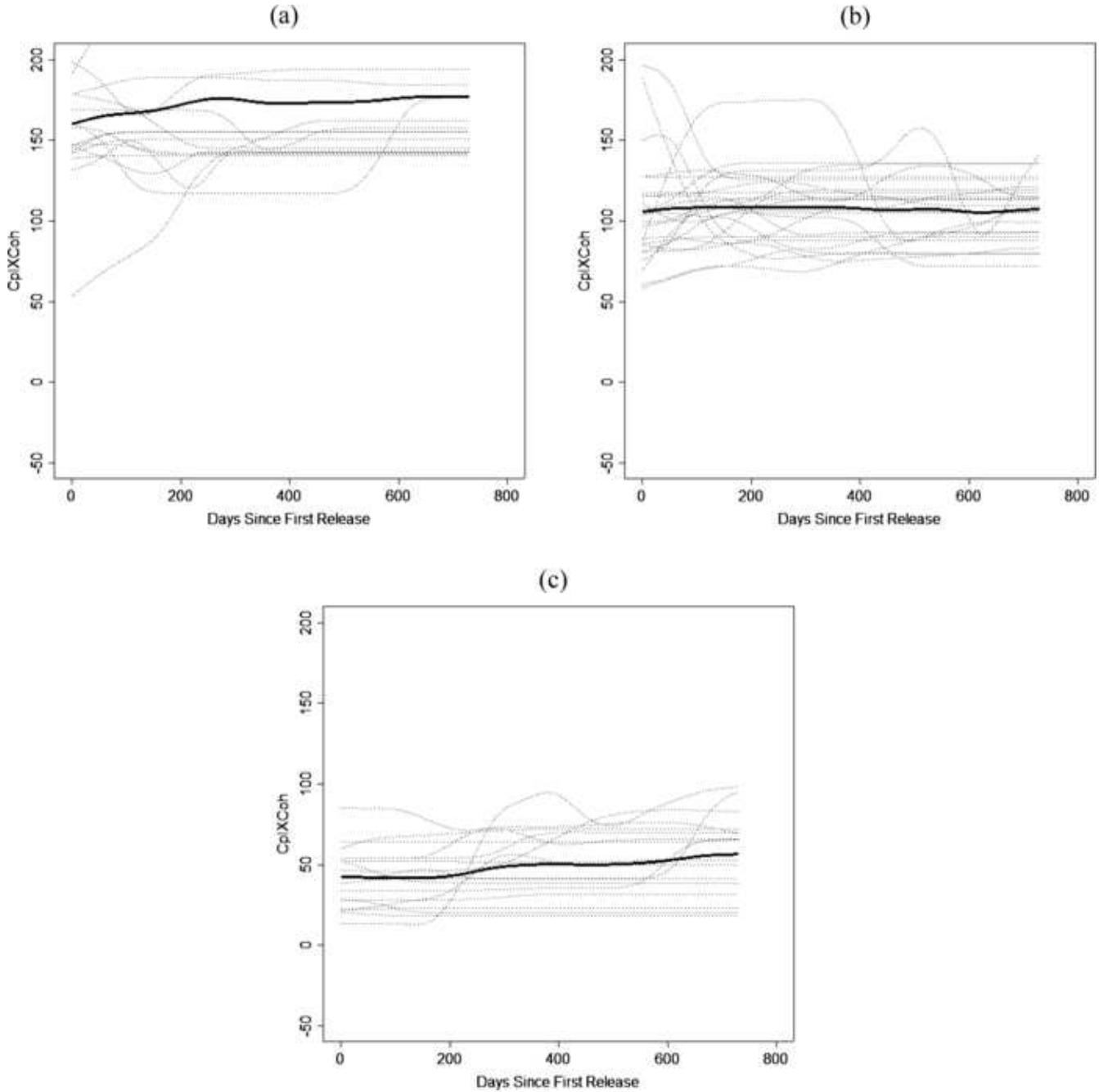

FIG. 4.   (a) *Cluster 1 of 3 (absolute values). (b) Cluster 2 of 3 (absolute values). (c) Cluster 3 of 3 (absolute values).*

The three-cluster solution [Figure 5(b)] produces a decreasing-complexity cluster very similar to that seen in the two-cluster solution, while the increasing-complexity cluster splits into two clusters with different patterns. One cluster of increasing-complexity projects is active over approximately the first 150 days and then flattens out. The second cluster shows an upward trend throughout the observation period, with the fastest rate of increasing complexity between approximately days 200 and 400. These projects appeared to have longer active periods than the projects in either of the other clusters.

Table 4
*Descriptive statistics for the four-cluster solution*

|  | Cluster 1— early decreasers | Cluster 2— early increasers | Cluster 3— midterm increasers | Cluster 4— midterm decreasers |
|---|---|---|---|---|
| Number of releases | 6.23 | 5.94 | 12.21 | 9.93 |
| Average release frequency (in days) | 40.93 | 36.01 | 80.60 | 71.42 |
| First release LOC | 4,893 | 6,576 | 5,052 | 11,579 |
| Last release LOC | 10,872 | 14,166 | 13,934 | 19,139 |
| First release CplxLCoh | 113.99 | 104.49 | 58.99 | 114.97 |
| Last release CplXLCoh | 111.31 | 140.21 | 99.78 | 95.49 |

Table 3
*Project cluster membership*

|  | Two clusters | Three clusters | Four clusters |
|---|---|---|---|
| Cluster 1—early decreasers | 25 | 24 | 13 |
| Cluster 2—early increasers | 34 | 20 | 18 |
| Cluster 3—midterm increasers |  | 15 | 14 |
| Cluster 4—midterm decreasers |  |  | 14 |

The four-cluster solution [Figure 5(c)] resulted in a similar split in the decreasing-complexity projects. One decreasing-complexity cluster shows a relatively small change over the first 200 days, then flattens out before increasing during the last 150 days. The other cluster shows a larger decrease in complexity over approximately the first 400 days. The five-cluster solution did not result in any new pattern becoming apparent; therefore we focus on the four-cluster solution. For the four-cluster solution, the average within-cluster distance was 6.04 and the average intercluster distance was 3.60.

The groupings identified in the four-cluster solution might be differentiated as *early decreasers* that experience decreasing complexity during the first several months of the observed period, appear relatively stable during the middle of the observed period and then start trending back upward toward the end of the two years (cluster 1); *early increasers* that experience their fastest rate of increasing complexity during the first several months of the project and appear relatively stable after day 150 (cluster 2); *midterm increasers* that experience their fastest rate of increasing complexity several months after the start of the project and then appear more stable (cluster 3); and *midterm decreasers* that continue a steady decrease through the midpoint of the period before stabilizing around day 400 (cluster 4).

The main differentiating characteristics of the projects across clusters appear to be the change in complexity (i.e., increased or decreased) and the period of time over which the projects were most active. To compare the four-cluster solution to the three- and five-cluster solutions and to confirm that apparent differences across clusters are representative of real differences across the projects and not an artifact of the smoothing splines, we conducted multiple analysis of variance (MANOVA). The MANOVA used the percentage change in CplXLCoh for each project and the active life of each project as dependent variables with the cluster assignment as the categorical independent variable. Percentage change in CplxLCoh was calculated for each project using the non-standardized data: (CplXLCoh of the last release – CplXLCoh of the first release)/CplXLCoh of the first release. Active life was calculated as the number of days between the first release and the last release of a project during the 730 days over which projects were observed.

For each of the three-, four- and five-cluster solutions, cluster assignment had a significant impact on both of the dependent variables ($p < 0.01$ in all cases). The three-cluster solution explained 17.6% and 17.3% of the variance in percentage change in CplXLCoh and active life, respectively. In the four-cluster solution these numbers were increased to 20.0% and 27.9%. In the five-cluster solution the explained variance decreased slightly to 19.9% and 27.5%. These results seem to be consistent with our conclusion above that the four-cluster solution appears to be of greatest interest. For the four-cluster solution, the multivariate $F$-test was significant ($F_{3,55} = 7.582$, $p < 0.001$) as were univariate tests for each dependent variable (for percentage change in CplXLCoh $F_{3,55} = 4.570$, $p < 0.01$, and for active life $F_{3,55} = 7.081$, $p < 0.001$). Differences across clusters were in the expected directions in all cases.

Given the relatively small sample size, and the limited amount of other information available in our dataset about the projects, it is difficult to identify what other characteristics may be associated with cluster membership. Prior work has suggested there may be correlations between size in LOC and complexity or between the release pattern of a project and complexity (Tan and Mookerjee, 2005). Based on these suggestions, we examined four characteristics of potential interest in an attempt to uncover such differences. These were the size of the first release of a project, the number of releases, the average release frequency and the percentage increase in size. We entered these variables into a MANOVA with the cluster assignment as the predictor variable. The multivariate $F$-test was significant ($F_{4,54} = 9.085$, $p < 0.001$), as were univariate tests for number of releases ($F_{3,55} = 3.564$, $p < 0.05$) and average release frequency ($F_{3,55} = 2.976$, $p < 0.05$). Univariate tests for size of the first release and percentage increase in size were not significant at $\leq 0.05$. Pairwise comparisons indicated that projects in the midterm increasers cluster had significantly more releases than the early increasers (difference $= 6.27$, $p < 0.01$) or the early decreasers (difference $= 5.984$, $p < 0.05$). The average release frequencies for these projects were also longer than the average release frequencies for the early increasers (difference $= 44.589$, $p < 0.05$) or the early decreasers (difference $= 39.665$, $p < 0.05$). The midterm decreasers cluster had a significantly longer average release frequency than the early increasers (difference $= 35.410$, $p < 0.05$). All other pairwise comparisons were insignificant.

**Sensitivity Analysis**

As discussed above, several parameters were chosen in conducting the analysis. To reduce the likelihood that the results reported here are artifacts of specific parameter values, we repeated the analysis varying these parameters. The number of knots was changed to 6 and to 26. The smoothing parameter was varied from the minimum to the maximum allowable value for the smoothing function that we used to generate the splines. We performed clustering using both the coefficients and the predicted values for the curves. In no case were the conclusions from the analysis affected. While minor variations in the shapes of the cluster mean curves were observed, in every case we observed four clusters with the same basic shapes (i.e., two clusters of increasing-complexity projects and two clusters of decreasing-complexity projects, as described).

In addition to varying the parameters, we applied different screening criteria to the data, screening out projects with fewer than three releases and including projects with any positive increase in LOC. Again, while minor variations in the shapes of the curves were observed, conclusions were unaffected.

# 6. DISCUSSION

Though this initial analysis has focused on a small sample of projects, results indicate FDA is a promising approach for studying the dynamics of evolution in software development. Identifying different patterns of change is a first step toward building a better understanding of why projects fall into one pattern versus another and what differences in outcomes may result from different patterns. This is of particular interest given that patterns of decreasing complexity are likely to result in desirable outcomes such as lower maintenance costs (Darcy, Kemerer et al., 2005).

Preliminary explorations into why some of the projects grouped into decreasing (complexity) clusters and others into increasing (complexity) clusters indicate that some of the intuitive and long-held beliefs about project complexity may not apply in some OSS development projects. In particular, neither the starting size nor the increase in size over the observed history of a project was significantly different across increasing-complexity or decreasing-complexity clusters. This is surprising given the widely held belief that increasing size is associated with increasing complexity (Chidamber, Darcy and Kemerer, 1998). Similarly, the pattern of differences across clusters in the average release frequency of projects did not distinguish between increasing- and decreasing-complexity clusters, but rather between clusters that had longer and shorter active periods. In the future research section below, we discuss other factors that may influence the assignment of a project to a cluster and the potential implications of cluster membership.

## Limitations

Aside from possible limitations on the generalizability of the particular patterns uncovered here (due to the screening criteria used to select projects), the decisions regarding how the data were analyzed have important implications. First, we focused on CplXLCoh as a measure of complexity because it has been shown in prior work (Darcy, Kemerer et al., 2005) to better represent overall structural complexity than either Cpl or LCoh alone. The individual measures, Cpl and LCoh, or other measures such as size or McCabe's Cyclomatic complexity may exhibit different patterns over time. Second, our approach to creating comparable data across projects by interpolating values of complexity for every day could underestimate the uncertainty regarding the levels of complexity and trajectory of changes in projects with very few releases. For example, the pointwise calculation of the 95% confidence interval using fitted values creates a relatively tight bound around some areas of the mean curve where there may be relatively few actual observations (i.e., releases). The approach taken in this analysis was to weight every project equally whereas for some projects we have many more data points than for others.

It is important to bear in mind that we cannot interpret the curves as representing what happens during the first two years of software development, but rather during the first two years of *public* development. That is, the results tell us about the nature of the first two years of evolution *after* a project is opened to the community, and that, as explained above, could happen at different stages of a project's lifecycle. One interesting avenue for future study is to investigate whether a reason projects may fall into the different clusters uncovered in this analysis is due to their having initially released software at different stages. For example, it could be that projects in clusters where complexity decreased released software with all of the intended functionality already included and the main effort of the community was then to "clean it up" and add minor extensions to an already well-defined architecture. In contrast, perhaps the projects in the increasing-complexity clusters released very early initial versions and thus the first two years of open development represent more significant additions of functionality. However, if size represents functionality, then the exploratory post hoc analysis showing no statistical difference between the percent change in size across clusters would cast some doubt on this possibility.

As noted above, projects were only observed for approximately two years. Based on the results it seems that this may be long enough to capture the dynamics of many projects because all cluster curves exhibited the most activity prior to approximately day 400. However, while our screening criteria for inclusion in the sample ensured projects did maintain a presence on Sourceforge during the entire two years of observation, it is possible that some projects may maintain such a presence while shifting development work to another location, which would be an alternative reason that development might appear to cease. It remains to be seen if the projects in the more dynamic upward trending cluster may later reach equilibrium or if they will continue to grow steadily in complexity over time.

## Future Research

There are several ways to build on the exploratory work presented here. Many of these revolve around refining the patterns uncovered, examining the extent to which they may be replicated in other samples (e.g., projects that use different programming languages) and exploring their antecedents and consequences. One potential antecedent to cluster membership mentioned above is the development stage of the project. Prior work on software development has indicated that larger development teams may produce more complex software (Banker, Davis and Slaughter, 1998); thus one fruitful avenue for future work on antecedents of cluster membership may be to examine the project team size. Because OSS projects often rely on voluntary labor, team size is an issue of particular relevance in the OSS setting and attracting larger teams has generally been viewed as an indicator of project success (Stewart and Gosain, 2006). Projects that experience decreasing complexity may be able to attract and retain more developers because new developers will be able to more quickly understand and modify the code. Combined with the finding in prior work that larger teams produce more complex code, this could explain the pattern in the early decreasers cluster in which complexity initially decreases (allowing for the attraction of more developers) but then trends upward toward the end of the observation period (once the team size has increased).

In addition to having different sizes, OSS projects have been observed to use different kinds of organizational structures (e.g., Apache has a voting system in place to advance contributors into positions of leadership whereas many other projects have less formalized processes), and prior work has suggested that the organization of software code may mirror the organization of the group that produces it (Mac-Cormack, Rusnak and Baldwin, 2004). Thus another potential antecedent to cluster membership may be the organizational structure of the team that develops the software.

A reason that uncovering antecedents to cluster membership may be important is that cluster membership may have implications for the future success of projects. For example, Yu, Schach et al. (2004) postulate that the currently high complexity of Linux may shortly bring it to a point where maintenance becomes extremely difficult. Generalizing from that work, it may be the case that projects in the increasing-complexity clusters may eventually be unable to add enhancements to keep up with changing user needs and thus lose popularity. Given that popularity has been considered a facet of OSS project success (Stewart, Ammeter and Maruping, 2006), such a trend would be undesirable for the longer-term health of an OSS project.

In addition to these questions for which there is a ready arsenal of analytical techniques to apply, there are other questions that may require more sophisticated statistical strategies. For example, because structural complexity has multiple components, we used a cross-term of two such components to capture an overall level of complexity (Cpl X LCoh). Since there is a trade-off between Cpl and LCoh, it would be interesting to be able to study their co-evolution—that is, do projects that manage to contain overall complexity (e.g., those in the downward sloping clusters) do so by alternately attending to LCoh and Cpl or by attending to both in a consistent manner over time? Similarly, to what extent do changes in size co-occur with changes in complexity? A means for analyzing patterns in two or more curves for each project simultaneously, such as a bivariate functional object, may be helpful in addressing such questions.

It may also be of interest to attempt to replicate or refine the findings from this study using different approaches. For example, an alternative analytical approach would be to use functional principal components analysis rather than cluster analysis. An alternative approach to examining evolution may be to analyze more granular data. For example, rather than focusing on official releases, one might track each change as it occurs through some form of version or configuration management system (such as the CVS mechanism on Sourceforge).

## Conclusion

The goal for this paper was to explore FDA as a means to uncover and characterize patterns of evolution in the complexity of OSS projects. FDA enables the examination of projects through the creation and manipulation of a functional data object for each project. This approach allows a richer exploration and comparison of projects than has been previously possible using prior approaches in the software engineering literature. This analysis has suggested insights into the evolution of complexity in open source projects, particularly with regard to the existence of multiple nuanced patterns. Our post hoc explorations have suggested that there is a substantial need for extensive additional research into the various patterns. These initial findings open extensive avenues for further research that we hope will ultimately provide practical guidance in managing software complexity.

## ACKNOWLEDGMENTS

We are grateful for the research assistance of Julie Inlow, Chang-Han Jong and Vincent Kan as well as the helpful guidance of Wolfgang Jank and Galit Shmueli. We also appreciate the constructive suggestions of participants at the First Interdisciplinary Symposium on Statistical Challenges in e-Commerce. This project was supported by the National Science Foundation Grant IIS-03-47376. Any opinions, findings and conclusions or recommendations expressed in this material are those of the authors and do not necessarily represent the views of the National Science Foundation. Data described in this paper, including a listing of projects studied, is available by contacting the first author.

## REFERENCES

Banker, R. D., Davis, G. B. and Slaughter, S. A. (1998). Software development practices, software complexity, and software maintenance performance: A field study. *Management Sci.* **44** 433–450.

Belady, L. A. and Lehman, M. M. (1976). A model of large program development. *IBM Systems J.* **15** 225–252.

CHIDAMBER, S. R., DARCY, D. P. and KEMERER, C. F. (1998). Managerial use of metrics for object-oriented software: An exploratory analysis. *IEEE Trans. Software Engineering* **24** 629–639.

DARCY, D. P., KEMERER, C. F., SLAUGHTER, S. A. and TOMAYKO, J. E. (2005). The structural complexity of software: An experimental test. *IEEE Trans. Software Engineering* **31** 982–995.

GORLA, N. and RAMAKRISHNAN, R. (1997). Effect of software structure attributes on software development productivity. *J. Systems and Software* **36** 191–199.

HASTIE, T., TIBSHIRANI, R. and FRIEDMAN, J. (2001). *The Elements of Statistical Learning*. Springer, New York. MR1851606

JANK, W. and SHMUELI, G. (2005). Profiling price dynamics in online auctions using curve clustering. Working paper RHS-06-004, Smith School of Business, Univ. Maryland. Available at ssrn.com/abstract=902893.

KAUFMANN, L. and ROUSSEEUW, P. J. (1987). Clustering by means of medoids. In *Statistical Analysis Based on the $L_1$-Norm and Related Methods* (Y. Dodge, ed.) 405–416. North-Holland, Amsterdam.

KEMERER, C. F. (1995). Software complexity and software maintenance: A survey of empirical research. *Annals of Software Engineering* **1** 1–22.

KEMERER, C. F. and SLAUGHTER, S. A. (1999). An empirical approach to studying software evolution. *IEEE Trans. Software Engineering* **25** 493–509.

MACCORMACK, A., RUSNAK, J. and BALDWIN, C. (2004). Exploring the structure of complex software designs: An empirical study of open source and proprietary code. Working paper 05-016, Harvard Business School.

PRAHALAD, C. K. and KRISHNAN, M. S. (1999). The new meaning of quality in the information age. *Harvard Business Review* Sept. 109–118.

RAMSAY, J. O. and SILVERMAN, B. W. (2002). *Applied Functional Data Analysis*: *Methods and Case Studies*. Springer, New York. MR1910407

SCACCHI, W. (2002). Understanding the requirements for developing open source software systems. *IEEE Proc. Software* **149** 24–39.

SHMUELI, G. and JANK, W. (2006). Modeling the dynamics of online auctions: A modern statistical approach. In *Economics, Information Systems and E-Commerce Research II*: *Advanced Empirical Methods* **1** (R. Kauffman and P. Tallon, eds.). Sharpe, Armonk, NY. To appear.

SMITH, T. (2002). Open source: Enterprise ready—with qualifiers. Available at www.linuxtoday.com/it_management/2002100101126NWBZ.

STEWART, K., AMMETER, A. and MARUPING, L. M. (2006). Impacts of license choice and organizational sponsorship on user interest and development activity in open source software projects. *Information Systems Research* **17** 126–144.

STEWART, K. and GOSAIN, S. (2006). The impact of ideology on effectiveness in open source software development teams. *Management Information Systems Quarterly* **30** 291–314.

TAN, Y. and MOOKERJEE, V. S. (2005). Comparing uniform and flexible policies for software maintenance and replacement. *IEEE Trans. Software Engineering* **31** 238–255.

YU, L., SCHACH, S. R., CHEN, K. and OFFUTT, J. (2004). Categorization of common coupling and its application to the maintainability of the Linux kernel. *IEEE Trans. Software Engineering* **30** 694–706.